\documentclass[hypertex,11pt]{article}
\usepackage{amsfonts,amsmath,amscd,amssymb,amsthm}
\usepackage{fullpage}
\usepackage[alphabetic,initials]{amsrefs}
\usepackage[all]{xy}
\usepackage{hyperref} 
\usepackage{mathrsfs}
\input xy%
\usepackage{fullpage}
\usepackage{amssymb}

\def\convf{\hbox{\space \raise-2mm\hbox{$\textstyle      \bigotimes \atop \scriptstyle \omega$} \space}}
\def\0{{\bar 0}}
\def\1{{\bar 1}}

\def\Z{{\mathbb Z}}

\def\Res{{\operatorname{Res}}}
\def\Ind{{\operatorname{Ind}}}

\def\ch{{\operatorname{ch}\:}}

\def\pwg1{{\operatorname{PWG}}}
\def\pwg{{\operatorname{pwg}}}

\def\span{{\operatorname{span}}}
\def\Tor{{\operatorname{Tor}}}

\def\Hom {{\operatorname{Hom}}}

\def\Mod {{\operatorname{Mod-}}}

\def\ann {{\operatorname{ann}}}

\newcommand{\ttk}{\mathtt{k}}
\newcommand{\tte}{\mathtt{e}}

\newcommand{\gL}{\Lambda}

\newcommand{\itema}{\item[{{\rm(a)}}]}
\newcommand{\itemb}{\item[{{\rm(b)}}]}

\newcommand{\noi}{\noindent}
\newcommand{\ga}{\alpha}
\newcommand{\gb}{\beta}
\newcommand{\gc}{\gamma}

\newcommand{\Gd}{\Delta}

\newcommand{\gd}{\delta}

\newcommand{\gs}{\sigma}

\newcommand{\gl}{\lambda}

\newcommand{\gr}{\rho}

\newcommand{\gep}{\epsilon}

\newcommand{\op}{\oplus}

\newcommand{\ot}{\otimes}

\newcommand{\fg}{\mathfrak{g}}\newcommand{\fgl}{\mathfrak{gl}}

\newcommand{\fh}{\mathfrak{h}}

\newcommand{\fm}{\mathfrak{m}}
\newcommand{\fn}{\mathfrak{n}}

\newcommand{\fc}{\mathfrak{c}}

\newcommand{\fk}{\mathfrak{k}}
\newcommand{\fp}{\mathfrak{p}}

\newcommand{\fl}{\mathfrak{l}}

\newcommand{\sfp}{\small{\mathfrak{p}}}

\newcommand{\sfg}{\small{\mathfrak{g}}}

\newcommand{\ff}{\footnote}
\newfont{\eufm}{eufm10 scaled\magstep1}

\newcommand{\ci}{\circ}

\newcommand{\cO}{\mathcal{O}}
\newcommand{\cC}{\mathcal{C}}

\newcommand{\cN}{\mathcal{N}}

\newcommand{\cF}{\mathcal{F}}

\newcommand{\cL}{\mathcal{L}}

\newcommand{\cS}{\mathcal{S}}
\newcommand{\cT}{\mathcal{T}}

\newcommand{\ey}{\end{eqnarray}}
\newcommand{\by}{\begin{eqnarray}}
\newcommand{\nn}{\nonumber}

\newcommand{\bco}{\begin{conjecture}}
\newcommand{\ba}{\begin{alg}}
\newcommand{\ea}{\end{alg}}
\newcommand{\eco}{\end{conjecture}}
\newcommand{\bpf}{\begin{proof}}
\newcommand{\epf}{\end{proof}}
\newcommand{\bt}{\begin{theorem}}

\newcommand{\et}{\end{theorem}}

\newcommand{\br}{\begin{rem}}
\newcommand{\er}{\end{rem}}
\newcommand{\brs}{\begin{rems}}
\newcommand{\ers}{\end{rems}}
\newcommand{\bi}{\begin{itemize}}
\newcommand{\ei}{\end{itemize}}
\newcommand{\bl}{\begin{lemma}}
\newcommand{\bsul}{\begin{sublemma}}
\newcommand{\esul}{\end{sublemma}}
\newcommand{\bp}{\begin{proposition}}
\newcommand{\be}{\begin{equation}}
\newcommand{\bc}{\begin{corollary}}
\newcommand{\bexs}{\begin{examples}}
\newcommand{\eexs}{\end{examples}}
\newcommand{\bexa}{\begin{example}}
\newcommand{\eexa}{\end{example}}
\newcommand{\bex}{\begin{exercise}}
\newcommand{\eex}{\end{exercise}}
\newcommand{\btab}{\begin{tab}}
\newcommand{\etab}{\end{tab}}
\newcommand{\el}{\end{lemma}}
\newcommand{\ep}{\end{proposition}}
\newcommand{\ee}{\end{equation}}
\newcommand{\ec}{\end{corollary}}
\newcommand{\Bc}{\begin{center}}
\newcommand{\Ec}{\end{center}}
\newcommand{\bh}{\begin{hyp}}
\newcommand{\eh}{\end{hyp}}
\newcommand{\bhs}{\begin{hyps}}
\newcommand{\ehs}{\end{hyps}}
\newcommand{\bd}{\begin{dfn}}
\newcommand{\ed}{\end{dfn}}

\begin{document}
\title{Table of Contents}

\newtheorem{thm}{Theorem}[section]
\newtheorem{hyp}[thm]{Hypothesis}
 \newtheorem{hyps}[thm]{Hypotheses}
  \newtheorem{rems}[thm]{Remarks}

\newtheorem{conjecture}[thm]{Conjecture}
\newtheorem{theorem}[thm]{Theorem}
\newtheorem{theorem a}[thm]{Theorem A}
\newtheorem{example}[thm]{Example}
\newtheorem{examples}[thm]{Examples}
\newtheorem{corollary}[thm]{Corollary}
\newtheorem{rem}[thm]{Remark}
\newtheorem{lemma}[thm]{Lemma}
\newtheorem{sublemma}[thm]{Sublemma}
\newtheorem{cor}[thm]{Corollary}
\newtheorem{proposition}[thm]{Proposition}
\newtheorem{exs}[thm]{Examples}
\newtheorem{ex}[thm]{Example}
\newtheorem{exercise}[thm]{Exercise}
\numberwithin{equation}{section}%
\setcounter{part}{0}
\newcommand{\drar}{\rightarrow}
\newcommand{\lra}{\longrightarrow}
\newcommand{\rra}{\longleftarrow}
\newcommand{\dra}{\Rightarrow}
\newcommand{\dla}{\Leftarrow}

\newtheorem{Thm}{Main Theorem}


\newtheorem*{thm*}{Theorem}
\newtheorem{lem}[thm]{Lemma}
\newtheorem*{lem*}{Lemma}
\newtheorem*{prop*}{Proposition}
\newtheorem*{cor*}{Corollary}
\newtheorem{dfn}[thm]{Definition}
\newtheorem*{defn*}{Definition}
\newtheorem{notadefn}[thm]{Notation and Definition}
\newtheorem*{notadefn*}{Notation and Definition}
\newtheorem{nota}[thm]{Notation}
\newtheorem*{nota*}{Notation}
\newtheorem{note}[thm]{Remark}
\newtheorem*{note*}{Remark}
\newtheorem*{notes*}{Remarks}
\newtheorem{hypo}[thm]{Hypothesis}
\newtheorem*{ex*}{Example}
\newtheorem{prob}[thm]{Problems}
\newtheorem{conj}[thm]{Conjecture}

\title{ Twisting Functors and Generalized Verma modules.}
\author{Ian M. Musson\ff{Research partly supported by  Simons Foundation grant 318264.} \\Department of Mathematical Sciences\\
University of Wisconsin-Milwaukee\\ email: {\tt
musson@uwm.edu}}
\maketitle
\begin{abstract}
Let $\fg$ be a reductive Lie algebra. We give a condition that ensures that the character of a generalized Verma module is well-behaved under a twisting functor. We show that a similar result holds for basic classical simple Lie superalgebras.  The result is used in \cite{M1001} to obtain a Jantzen sum formula for  certain  
highest weight modules over type A Lie superalgebras. 
\end{abstract}
\section{Introduction and Notation.}
\subsection{Introduction.} Let $\ttk$ be an algebraically closed field of characteristic zero, and 
$\fg$ a reductive Lie algebra over $\ttk$. Our current understanding of the BGG category $\cO$ of $\fg$-modules owes much to  several well known endofunctors of the category $\cO$ and its blocks, see \cite{BGG}, \cite{H2}, \cite{Mz}.  
An important class of functors are the twisting functors introduced by
Arkhipov \cite{A} and studied further by Andersen in collaborations with Lauritzen \cite{AL}, and Stroppel \cite{AS}. 
Twisting functors for Lie superalgebras are studied in \cite{CM}.  
We postpone some definitions until later, but remark  here that there is a twisting functor $T_w$ for any Weyl group element $w$, and for any Verma module $M(\gl)$, the modules $T_wM(\gl)$ and $M(w\cdot \gl)$ have the same character.  
\\ \\
Suppose  $\fp$ is a parabolic subalgebra of $\fg$ with Levi factor $\fl$. 
Let  $\cO_\fl$  be the 
analog of the category of $\cO$ for $\fl$.  
Consider the functor $\cF:\cO_\fl\lra \cO$ 
defined on an object $L$ as follows. 
First make $L$ into a $\fp$-module by allowing the radical $\fm^+$ of $\fp$ to act trivially, and then let 
$\cF L =\Ind_{\fp}^{\fg} L$  be the induced $\fg$-module. 
If $L$ is finite dimensional simple, we call $\Ind_{\fp}^{\fg} L,$ a generalized Verma module. 
An analog of the above mentioned result on twisted Verma modules does not always hold for generalized Verma modules.
However if $w\ga$ is a positive root  of $\fg$ for every positive root $\ga$ of  $\fl$, there is a suitable analog for the character of $T_w\Ind_{\fp}^{\fg} L$.  The result is most clearly expressed using partition functions.  
Suppose that $X$ is any set of positive roots, and define
$p_X = \prod_{\ga \in \Gd^+ \backslash X} (1 -\tte^{-\ga})^{-1}$. If $\fl$ and $\fp$ are as above and $X$ is the set of positive roots of $\fl$, then the character of $\cF L=\Ind_{\fp}^{\fg} L$ is given by 

\[\ch\Ind_{\fp}^{\fg} L = \ch L\;p_X .\]
Recall that the Verma module $M(\gl)$ with highest weight $\gl$ has character $\tte^\gl p$ where $p$ is  the Kostant  partition function. 
We can think of $p_X$ as the generating function for the number of partitions with support disjoint from $X.$ When $X$ is empty we have $p_X=p.$
Under the above assumption on $w$ we show that $\ch T_w \cF L = w\cdot \ch L\;p_{wX}$. In particular if $L$ is one-dimensional with weight 
$\gl\in \fh^*$ such that $\gl([\fl,\fl]\cap \fh)=0$, then $\ch T_w \cF L = \tte^{w\cdot \gl}p_{wX}$. 
\\ \\   
Furthermore the result extends easily to the case where $\fg$ is a classical simple Lie superalgebra, with a simple change in the definition of $p_X$ to accommodate odd roots. In \cite{M1001}
we showed in the contragredient case, that if $X$ is an isotropic set of (necessarily odd) roots, and $\gl\in\fh^*$ satisfies $(\gl+\gr,X) =0$, then there is a highest weight module $M^X(\gl)$ with highest weight $\gl$ such that $\ch M^X(\gl) = \tte^{ \gl}p_{X}$.  In type A, we obtained a Jantzen sum formula in the Grothendieck group $K(\cO)$ of the category $\cO$   for these modules.  The sum of 
terms in the Jantzen filtration is expressed as a linear combination, {\it with positive coefficients,} of other modules corresponding to a factorization of the \v Sapovalov determinant.  Among the modules that arise are some with character $\tte^\mu p_Y$ where $Y$ is a set of roots containing $X$ and 
the positive roots of a certain $\fgl(2|2)$-subalgebra of $\fg$.  
It does not seem easy to construct such modules using the methods from \cite{M1001}.  However they can be constructed easily by applying twisting functors to generalized Vermas. 
\\ \\
More generally, suppose that $\fk,\fl$ are subalgebras of $\fg$, such that 
$\fk$ is a contragredient Lie superalgebra or a reductive Lie algebra, 
and for some $r,s, t$
\be\label{fox}\fl \cong \fk \op \fgl(1|1)^r\op\fgl(1|0)^{s}\op \fgl(0|1)^{t}.\ee
We assume that $\fl$ contains a Cartan subalgebra $\fh$ of $\fg$. Then if 
$Y=\Gd^+(\fl)$, and $\gl\in \fh^*$ defines a one-dimensional representation $\ttk_\gl$ of $\fl$, 
the module $M^{Y}(\gl):=\Ind_{\fp}^{\fg} \ttk_\gl$ has character $\tte^\gl p_Y.$ 
Note that the set $V$ of positive roots in the $r$ copies of $\fgl(1|1)$ in \eqref{fox} is isotropic, and $V$ is orthogonal to all roots of $\fk$.  
Thus $X$ could consist of the union of $V$ and an isotropic set of roots of $\fk$. If $\fk=\fgl(2|2)$, 
and $|X|=r+2$, twisting  the induced modules $M^{Y}(\gl)$ gives the missing modules
referred to in the previous paragraph.  
Since the result seems to be new even for the semisimple case, and it holds in  greater generality 
than is needed for \cite{M1001},  it seems worth recording independently.

\subsection{Notation.} 
We collect the main notation that will be used in this paper. 
We assume that $\fg$ is a reductive Lie algebra or a contragredient Lie superalgebra.
If $\fg$ is a type A Lie superalgebra, we assume that $\fg=\fgl(m|n)$.  
Fix a Cartan subalgebra
$\fh$ of $\fg$, and  let 
$\fp$ be a parabolic subalgebra of $\fg$, that is a subalgebra containing an $\fh$-stable Borel subalgebra. 
There are subalgebras $\fm^\pm$, and a Levi factor $\fl$ such that

\be \label{kid} \fg= \fm^-\oplus \fl \oplus \fm^+\ee
 and $\fp= \fl \oplus \fm^+.$ We assume that the positive roots of $\fl$ are positive in $\fg$. 
Set $\fh_\fl= [\fl,\fl] \cap \fh$. For any $\fh$ stable, $\Z_2$-graded subalgebra $\fc$ of $\fg$, 
let $\;\Gd_0^+(\fc)$ and $\;\Gd_1^+(\fc)$ be the set of even and odd positive roots of $\fc$ respectively. 
We define $\;\Gd^+(\fc)=\Gd_0^+(\fc)\cup \Gd_1^+(\fc)$, and $\Gd^+=\Gd^+(\fg), 
\Gd_0^+=\Gd_0^+(\fg), \Gd_1^+=\Gd_1^+(\fg)$. 
Set $$\gr_0 = \frac{1}{2}\sum_{\ga\in \Gd^+_0}\ga, \quad \gr_1 = \frac{1}{2}\sum_{\ga \in \Gd^+_1} \ga$$  
and  $\gr=\gr_0-\gr_1$. 
For any root $\ga$, we can choose a root vector $e_\ga\in\fg^\ga$ such that 
$\fg^\ga =  \ttk e_\ga$.  For an element $w$ of the Weyl group  $W$, set  
$N(w) = \{ \alpha \in \Delta_0^+ | w \alpha < 0 \}.$
We define translated actions of $W$ on $\mathfrak{h}^*$
by

\be\label{nag} w \cdot \lambda = w(\lambda + \rho) - \rho, \quad \quad
 w \circ \lambda = w(\lambda + \rho_{0}) - \rho_{0}. \ee
If $X$ is a set of  positive  roots, let $X_0$ (resp. $X_1$) be the set of even (resp. odd) roots contained in $X$.  Then define

\[r_X = \prod_{\ga\in X_0}1- \tte^{-\ga}\; \mbox{ resp. } \;s_X= \prod_{\ga\in \Gd_1^+ \backslash X_1}1+ \tte^{-\ga}. 
\] 
Clearly for $w\in W$,
\be \label{rex} wr_X = r_{wX} \mbox{ and } ws_X = s_{wX}.\ee
We also set $r = r_{\Gd^+_0},$ and $p_X = r_Xs_Xp$. We remark that in the Lie algebra case, $r=s=t=0$ in \eqref{fox}, and $s_X = 1$. 
If $X=\Gd^+(\fl)$, then since $U(\fm^-) $ has  character given by 
\[\ch U(\fm^-) = \prod_{\ga\in \Gd^+_0(\fm^-)}(1- \tte^{-\ga})^{-1}\prod_{\ga\in \Gd_1^+(\fm^-)}1+ \tte^{-\ga}= p_X,\] 
we have 
\by\label{hog} \ch \Ind^\fg_\fp\; L = \ch L\;p_X.\ey 
If $L$ is one dimensional with weight $\gl$, then $\gl(\fh_\fl)=0$.  
In this case, we denote 
$L$ by $\ttk_\gl$, and for $h\in \fh$,
 we have  $(h-\gl(h))\ttk_\gl =0$. 
Then if $X=\Gd^+(\fl)$ we have 
\[\ch \Ind_{\fp}^{\fg} \;\ttk_\gl =\tte^{\gl}p_X.\]

\noi Let $\gL$ be the lattice of functions on $\fh^*$ such that $(\gs,\ga^\vee) \in \Z$ for all simple roots $\ga$.     
As in \cite{H} 22.5, we use the group ring $\Z[\gL]$ with $\Z$-basis the symbols $\tte^\gs$ with $\gs \in \gL$. 
The circle action of $w\in W$ on $\Z[\gL]$ is defined by  
\be \label{owl} w\ci \tte^\gs = \tte^{w\ci \gs}.\ee
(There is a similar dot action defined using $\gr$ in place of $\gr_0$). 
Using \eqref{owl} when $\gs \in \fh^*$, $w$ induces a map $\Z[\gs+\gL] \lra \Z[w\ci\gs+\gL]$. 
 We note that the circle action is not an action by algebra automorphisms.  Instead we have for $a, b \in \Z[\gL]$ and $w\in W$,
\be \label{don} w\ci  (ab) = (wa)(w\ci  b).\ee  
\noi With this notation, we can now state the main result on Lie algebras. 
\bt \label{key} Suppose $\fg$ is reductive and set $X=\Gd^+(\fl)$. 
If $w\in W$ is such that $N(w)\cap X = \emptyset$, the module $T_w \Ind_\fp^\fg L$ has character $(w\ci \ch L)p_{wX} $.  
In particular $\ch T_w \Ind_\fp^\fg \ttk_\gl=\tte^{w\ci \gl}p_{wX}$.\et
\noi This result extends to classical simple Lie superalgebras. In particular this gives a  new proof of a result of Coulembier and Mazorchuk, \cite{CM}  Lemma 5.5 about twisting Verma modules.\\ \\
I thank Kevin Coulembier for some helpful correspondence.  
\section{Reductive Lie algebras.}
Throughout this section we assume that $\fg$ is reductive, and set $U = U(\fg)$. 
 \subsection{Basics on Twisting Functors.} \label{cat} 
\noi 
The twisting functor on the category $\Mod U$ of $U$-modules
is defined  as follows: Let \; $\fn_w = \fn^-\cap w^{-1}(\fn^+)$, and let $N_w = U(\fn_w)$. 
We make  $\fg$ into  a $\Z$-graded Lie superalgebra $\fg =\bigoplus_{i\in \Z} \fg(i)$ such that  
$\fg(0) =\fh$ and $\fg(\pm 1) =\op \fg_{\pm\ga}$, 
where the sum runs over all simple roots $\ga$. This  grading restricts to a grading on  $N_w$. Let 
$(N_w^\star)_i = \Hom_\ttk((N_w)_{ -i},\ttk)$.  Then 
$N_w^\star =\bigoplus_{ i\in \Z} (N_w^\star)_i$
is the graded dual of $N_w$. 
The $U$-bimodule $\cS_w$ is defined as
$\cS_w = U\ot_{N_w} N_w^\star$.  
For a proof that $\cS_w$ is a $U$-bimodule   see \cite{A}. As a {\it right} $U$-module we have, see \cite{AL} 6.1,

\be \label{doe} \cS_w \cong N_w^\star \ot_{N_w} U.\ee
The {\it twisting functor} $T_w:\Mod U \lra \Mod U$ corresponding to $w \in  W$ is
defined by $T_w (?) = \cS_w \ot_U (?)$. The functor $T_w$ restricts to an endofunctor on $\cO$. 

\bl \label{yak} The functor $T_w$ is right exact, and has left derived functors $\cL_iT_w$ given by   
\be\label{dog} \cL_iT_w(?) = \Tor_i^{N_w}(N_w^\star, ?).\ee \el
\noi {\it Proof} 
This follows since by the definition and \eqref{doe}, 
$$ \quad \quad\quad\quad\quad\quad  T_w (?) = \cS_w \ot_U (?) = N^\star_w  \ot_{N_w} U\ot_U(?) = N^\star_w  \ot_{N_w} (?).  \quad \quad \quad \quad \quad \hfill \Box$$
\noi We end this section with a remark that is the basis for the definition of twisting functors in the super case. When $\fg$ is reductive and $w=s_\ga$ is a simple reflection, $T_w$ has the following easy description. 
Let $U_s$ denote the localization of $U$ at the negative root vector $e_{-\ga}$.
Then $U_{(s)} = U_s/U$ is a $U$-$U$-bimodule, and $T_s M \cong  U_{(s)}\ot_U M$.  
There is an inner automorphism $\phi=\phi_\ga$ of $\fg$ such that $\phi(\fg_\gb) = \fg_{s_\ga(\gb)}$ for all roots 
$\ga$, and $\phi(\fh) = \fh$. 
The action of $\fg$ on   $U_{(s)}$ is twisted by $\phi$, i.e. the action 
is given by $(x,u)\lra \phi(x)u$. 
It was shown in \cite{AL}, Remark 6.1 ii) that  

\be \label{paw} T_{ws} \cong T_wT_s \mbox{ if }  ws>w \mbox{ and } s \mbox{ is a simple reflection}.\ee 

\subsection{Twisting Generalized Verma Modules.}
It was shown in \cite{AL}, Proposition 6.1 ii), 
that for any  
Verma module $M(\gl)$ the modules $T_w M(\gl)$ and $M(w\circ \gl)$ have the same character.
Equivalently by  \cite{J1} Satz 1.11, we have in the Grothendieck group $K(\cO)$ of the category $\cO$ that

\be \label{eel} [T_w M(\gl)] =[M(w\circ \gl)].\ee 
Our goal is to obtain an analog of this result for generalized Verma modules.

\bl \label{ox} Suppose $w\in W$ and that the sequence of $\fg$-modules 

\[ 0 \lra Q_1 \lra Q_2 \lra Q_3 \lra 0,\]
is exact. Then provided that $Q=Q_3$ is a free $($or even flat$)$ $N_w$-module, the sequence 

\[ 0 \lra T_w Q_1 \lra T_w Q_2 \lra T_w Q_3 \lra 0,\]
is also exact. This holds for example if $Q$ is induced from a $\fp$-module. 
\el

\bpf The twisting functor $T_w$  is right exact, so it suffices to show that its left derived functor $\cL_1 T_w$ satisfies 
$\cL_1 T_w Q=0$. This follows from Lemma \ref{yak}.
\epf

\bl \label{pet} Suppose 
\[ 0 {\longrightarrow} M_{n} \stackrel{f_{n-1}}{\longrightarrow} M_{n-1} \lra \cdots \stackrel{f_1}{\longrightarrow} M_{1} \stackrel{f_0}{\longrightarrow} M_0 \lra M\lra  0,  \] is a long exact sequence such that for 
all $j\ge 0$ and $i\ge1$, $\cL_iT_w M_j =\cL_iT_w M =0$. Then the sequence 

\[ 0 {\longrightarrow} T_wM_{n} \stackrel{T_wf_{n-1}}{\longrightarrow} T_wM_{n-1} \lra \cdots \stackrel{T_wf_1}{\longrightarrow} T_wM_{1} \stackrel{T_wf_0}{\longrightarrow} T_wM_0 \lra T_wM\lra  0,  \] 
is also exact.\el
\bpf  Let $K_j = \ker f_j$. 
Then from the exact sequence $0\drar K_0 \lra M_0\lra M\drar 0$, and the long exact sequence for Tor, we see that 
$0\drar T_wK_0 \drar T_wM_0\drar T_wM\drar 0$ 
is exact and $\;\cL_iT_w K_0 = 0$ for all $i\ge 1$.  
The same reasoning applied inductively to the exact sequences $\;0\drar K_{j} \lra M_j\lra K_{j-1}\drar 0$ for 
 $j\ge 1$ shows that $0\lra T_wK_j \lra T_wM_j\lra T_wK_{j-1} \lra 0$ 
is exact and $\cL_iT_w K_j = 0$ for all $i, j$. 
Assembling the short exact sequences involving $T_w$ we get the result. 
\epf
\noi Let $\cC$ be the full subcategory of $\cO$ consisting of $\fg$-modules that have finite resolutions by (direct sums of) Verma modules, and let $\cC(\fl)$ 
be the analogous  category of 
$\fl$-modules. 
 The Verma module for $\fl$ with highest weight $\gl\in \fh^*$ will be denoted by $M_\fl(\gl)$. 
By associativity of the tensor product, $\cF M_\fl(\gl)\cong M(\gl)$, so the exact functor $\cF$ takes $\cC(\fl)$ to $\cC$. 

\bl \label{lox}If $L$ is a finite dimensional $\fl$-module, then $\cF L$ is an object  in $\cC$. \el
\bpf Begin with the BGG-resolution of $L$ as an $\fl$-module and apply $\cF$ to get a resolution of $\cF L$. 
As noted above, Vermas induce to Vermas.\epf 
\bp \label{bat} Suppose $N(w)\cap\Gd^+(\fl) = \emptyset$, then
\bi \itema $\fn_w$ is subalgebra of $\fm^-$.
\itemb If $L$ is a $\fl$-module and $M =\cF L$, then $\cL_i T_wM=0$ for $i\ge1$.
\ei \ep
\bpf By definition
\[\fn_w = \span \{e_{-\ga}|\ga \in N(w)\},\]
and \[\fm^- = \span \{e_{-\ga}|\ga \in \Gd^+(\fg)\backslash \Gd^+(\fl)\}.\]
Hence (a) follows from the hypothesis. With $L$ as in (b), $M=\Ind_\fp^\fg L= U(\fm^-)\ot L$ is a  free
$U(\fm^-)$-module. By the PBW Theorem and (a), $U(\fm^-)$ is a  free $N_w=U(\fn_w)$-module. Thus $M$ is free and hence 
flat over $N_w$, and (b) follows from \eqref{dog}.
\epf\noi 
\bl  \label{coq} Suppose $M\in \cC$, $w\in W$ and 
$\cL_i T_w M = 0$ for all $i\ge1$. Then
\bi \itema 
$\ch M = \sum_{\mu} b_\mu \ch M(\mu)$ implies 
$\ch T_w M = \sum_{\mu} b_\mu \ch M(w\ci \mu)$. 
\itemb If $\ch M = ap$ for $a\in \Z[\gL]$, then $\ch T_w M = (w\ci  a)p.$

\ei
\el \bpf \noi Suppose  that 

 \[ 0 {\longrightarrow} M_n {\longrightarrow} M_{n-1} \lra \cdots {\longrightarrow} M_{1} {\longrightarrow} M_0 \lra M\lra  0,  \] is a resolution of $M$ by direct sums of Vermas,  and suppose the multiplicity of the Verma $M(\mu)$ in $M_i$ is $a_{i,\mu}$, then clearly 

\by 
\ch M &=& \sum_{i=0}^n (-1)^i \ch M_i  = \sum_{i=0}^n \sum_{\mu}(-1)^i a_{i,\mu} \ch M(\mu)\nn\\
&=&\sum_{\mu} c_\mu \ch M(\mu)\nn,
\ey 
where $c_\mu =\sum_{i=0}^n (-1)^i a_{i,\mu}$.  
But since the classes of Verma modules form a  $\Z$-basis for $K(\cO)$, we have $c_\mu = b_\mu$.  
Now the Verma module $M(\mu)$ satisfies $\cL_iT_w M(\mu)=0$ 
for $i\ge1$, so $\cL_iT_w M_j =0$ if $j\ge 0$ and $i\ge1$. 
Since by assumption $\cL_i T_w M = 0,$  for all $i\ge1$,
Lemma \ref{pet} gives a resolution 

\[ 0 {\longrightarrow} T_w M_n {\longrightarrow} T_w M_{n-1} \lra \cdots {\longrightarrow} T_w M_{1} {\longrightarrow} T_w M_0 \lra T_w M\lra  0,  \] 
and  by \eqref{eel} this implies that  	

\[\ch T_w M = \sum_{i=0}^n \sum_{\mu} (-1)^i a_{i,\mu} \ch M(w\ci \mu)=\sum_{\mu} b_\mu \ch M(w\ci \mu).\]
This proves (a), and (b) follows by writing $\ch M = \sum_{\mu} b_\mu \ch M(\mu)$
and using (a). \epf

\noi {\it Proof of Theorem \ref{key}.} If $L$ is a finite dimensional $\fl$-module then $M= \Ind^{\sfg}_{\sfp} L\in \cC$ by Lemma \ref{lox}.
Also  $\cL_i T_wM=0$ for $i\ge1$ by Proposition \ref{bat}. 
Therefore from  \eqref{hog}, we have 

\[
\ch M=p_X \ch L=p r_{X}\ch L.\nn
\]
We  use Lemma \ref{coq} for the first equality below,
and then use \eqref{rex} and \eqref{don} to obtain
\by
\ch T_w \Ind^\fg_\fp  L =p{w\ci (r_X \ch L)}=p r_{wX}(w\ci \ch L),\nn
\ey
and this easily yields the result.\hfill  $\Box$
\section{Lie Superalgebras.}
From now on we fix a contragredient Lie superalgebra $\fg$, and  apply the results of the previous Subsection
to the reductive algebra $\fg_0$.  Recall Equations \eqref{fox} and \eqref{kid}. In this Section $p=1/r$, where $r$ is defined just after \eqref{rex}. If $X= \Gd^+(\fl)$ and $\gl(\fh_\fl)=0$, we set $M^X(\gl)= \Ind_{\fp}^{\fg} \;\ttk_\gl.$  In order to quote \eqref{hog} we need to modify the notation (since now $s_X\neq 1$).  Thus
if $L$  is a $\fp_0$-module, we have
\by\label{oie} \ch \Ind^{\fg_0}_{\fp_0}\; L = \ch L\;r_Xp.\ey
In \cite{CM} the twisting functor $T_s$   for $\fg$ is defined as follows.  We assume that $s$ is a  reflection corresponding to the simple non-isotropic root $\ga$. Then as before denote the localization $U(\fg)$   at $e_{-\ga}$ by $U_s$ and set $U_{(s)} = U_s/U$. 
Since the automorphism $\phi$ from Subsection \ref{cat} is inner, it extends to $\fg$ and then 
to $U(\fg)$. Then  the action of  $U$ is twisted by $\phi$.  
Moreover if we use \eqref{paw} to define an action 
of the free group generated by symbols $\cT_s$ for $s$ a simple reflection, then it is shown in \cite{CM} Lemma 5.3
that the  braid relations are satisfied. 
\ff{Braid relations for twisting functors are also discussed in \cite{KM} and \cite{P}.} Thus by choosing a shortest length expression  for  $w\in W$, \eqref{paw} yields a well-defined twisting functor (which we also denote by $T_w$).
Also the  
restriction functor $\Res^\fg_{\fg_0}$  intertwines $T_w$, that is we have an isomorphism of functors, see \cite{CM}	Lemma 5.1

\be \label{cow}
\Res^\fg_{\fg_0} \circ T_w  \cong T_w \circ \Res^\fg_{\fg_0}.
\ee
There is a similar result for the induction functor 
$ \Ind^\fg_{\fg_0}$, but we will not need it. 

\bl Suppose $E$ is a finite dimensional simple $\fl$-module, and make $E$ into a $\fp$-module by allowing $\fm^+$ to act trivially.  
\bi 
\itema
If $M= \Ind^{\sfg}_{\sfp} E$
there is a finite chain of $U(\fg_0)$-submodules

\be \label{sow} {M} = M_s \supset M_{s-1}\supset \dots \supset M_1 \supset M_0 = 0,\ee
and $\fl_0$-modules $E_i$ 
with $M_i/M_{i-1} \cong \Ind^{\sfg_0}_{\sfp_0} E_i$ and $\dim E_i <\infty$ for
$1 \leq i \leq s.$ 
\itemb Denote the exterior algebra on $\fm^-_1$  by $\gL(\fm^-_1)$. Then we have
\be \label{ant} \sum_{i=1}^s \ch E_i = \ch \gL(\fm^-_1)\ot E= s_X \ch E.\ee. 
\ei
\el

\bpf 
Note that $\fl_0$ acts on $\gL(\fm^-_1)$ via the adjoint action, and  we have as a $\fl_0$-module, (compare \cite{M} Corollary 6.4.5),

\by M = U(\fg)\ot_{U(\fp)} E &=& U(\fm^-)\ot E\nn\\
 &=& U(\fm^-_0) \ot \gL(\fm^-_1)\ot E.\nn\ey 
Note that $N=\gL(\fm^-_1)\ot E$ is a $\fp_0$-module, and $M\cong \Ind_{\fp_0}^{\fg_0} N$ as a $\fg_0$-module.   
Write $N$  as a direct sum of finite dimensional simple $\fl_0$-modules.
Say 

\be \label{gnu}N=\gL(\fm^-_1)\ot E = \bigoplus_{i=1}^s E_i \ee
with $E_i$ simple. Since $M$ is an object of $\cO$, 
$\fn_0^+$ acts nilpotently on $M$, so some power of $\fm_0^+$ in $U(\fm_0^+)$ annihilates $N$. 
To construct the chain \eqref{sow} we construct a similar chain of $U(\fp_0)$-submodules

\be \label{oie} {N} = N_s \supset N_{s-1}\supset \dots \supset N_1 \supset N_0 = 0\ee
such that $N_i/N_{i-1} \cong  E_i$ for $1 \leq i \leq s.$ Suppose we have constructed $N_i$ and set 
$M_i=\Ind^{\sfg_0}_{\sfp_0} N_i$ and $\cN_i= (M_i+N)/M_i$.  Then by the remark following \eqref{gnu}

\[\ann_{\cN_i} \fm_0^+ = \{x\in \cN_i|\fm_0^+ x = 0\}\] is non-zero.
Furthermore it is a $\fp_0$-submodule of $\cN_i$ killed by $\fm_0^+$, and so a module for $\fl_0$.  
By renumbering the $E_j$ in Equation \eqref{gnu} we can assume that $N_{i+1}/N_i\cong E_{i+1}$ is a simple submodule of $\ann_{\cN_i} \fm_0^+$. This completes the construction of the sequences \eqref{sow} and \eqref{oie}. 
Finally applying the induction functor  $\Ind^{\sfg_0}_{\sfp_0} $ to the exact sequence

\be  0\lra N_{i-1} \lra N_{i} \lra  E_i \lra 0,\nn\ee we obtain

\be \label{bee} 0\lra M_{i-1} \lra M_{i} \lra \Ind^{\sfg_0}_{\sfp_0} E_i \lra 0.\ee
Hence  (a) follows and (b) holds by \eqref{gnu}.
\epf  
\noi 

\bt \label{zoo} Assume $X= \Gd^+(\fl)$ and $w\in W$ is such that $N(w)\cap X  = \emptyset$. Then if $M= \Ind^{\sfg}_{\sfp} E$, we have 
\be \ch T_wM ={r_{wX} s_{wX}p(w\ci \ch E)} = p_{wX}(w\ci \ch E).\ee \et
\bpf
By \eqref{cow} the character of $T_w M$ is the same as its character when regarded as a $\fg_0$-module.  
\noi By Lemma \ref{ox}, the functor $T_w$ is exact on the sequence \eqref{bee}, so 
$T_w M_{i}/T_w M_{i-1} \cong T_w\Ind^{\sfg_0}_{\sfp_0}  E_{i}.$
These remarks justify the first two equalities below, and for the third we use \eqref{gnu}.  Thus we have 

\by \ch T_wM &=& \sum_{i=1}^s \ch T_wM_i/T_wM_{i-1} \nn \\
 &=& \sum_{i=1}^s  \ch T_w \Ind^{\sfg_0}_{\sfp_0} E_i\nn\\
&=& \ch T_w \Ind^{\sfg_0}_{\sfp_0} (\gL(\fm^-_1)\ot E),\nn
\ey 
but by \eqref{oie} and \eqref{ant},

\be \label{pup} \ch \Ind^{\sfg_0}_{\sfp_0} \gL(\fm^-_1)\ot E=
{r_{X}s_X p}\;\ch E .\ee
Hence using Lemma \ref{coq}, we obtain the result.
\epf 
\noi We specialize to the case where $E=\ttk_\gl$ is one dimensional.  
Then we give  an expression for $\ch T_w M^{X}(\gl)$ using the dot 
action rather than the circle action from Theorem \ref{zoo}. 
To do this we need some more notation. 
Let $\Gamma$ be the set of sums of distinct odd positive roots. 
By \cite{M0} Lemma 2.3, $W$ acts on $\Gamma$ by

\be\label{rat} w * \gamma = \rho_{1} - w(\rho_{1} - \gamma), \ee
for $w \in W$ and $\gamma \in \Gamma$.
This action is related to those in \eqref{nag} by 

\be \label{pig} w \circ (\lambda -\gamma) = w \cdot \lambda - w * \gamma.\ee
Let $Z$ be the set of odd roots in $X$ and for $\gc \in \Gamma,$ define $K_Z(\gc)$ by 

\be    \prod_{\alpha \in\Delta^{+}_{1} \backslash Z } (1 + \tte^{-\alpha}) = \sum_{\gamma\in \Gamma} K_Z(\gamma)\tte^{- \gamma} =s_X. 
\ee 
and let $\Gamma_Z = \{\gamma
\in \Gamma| K_Z(\gamma) > 0 \}.$
\noi 
Then  
\eqref{pup} yields
\by \label{cod} \ch \Ind^{\sfg_0}_{\sfp_0} \gL(\fm^-_1)\ttk_\gl  &=&{p_{X}}\tte^\gl 
\nn \\ &=& {r_{X}}p\sum_{\gc \in \Gamma_Z} K_Z(\gc)\tte^{\lambda -\gamma}.\ey


\bt \label{hen} Suppose $w\in W$ is such that $N(w)\cap\Gd^+(\fl) = \emptyset$.
Then $$\ch T_w M^{X}(\gl) =p_{wX}\tte^{w\cdot \gl}.$$\et
\noi First we isolate a key step in the proof.
\bl \label{pug} We have
$$\sum_{\gamma \in\Gamma_Z} K_Z(\gamma)\tte^{- w*\gamma} = \sum_{\gamma \in\Gamma_{wZ}} K_{wZ}(\gamma)\tte^{- \gamma}.  $$ 
\el

\bpf 
Note that

\be\label{asp} \tte^{\rho_{1}} \prod_{\alpha \in Z}(1 + \tte^{-\alpha})\sum_{\gamma \in \Gamma_Z}
K_Z(\gamma)\tte^{ - \gamma} = \prod_{\alpha \in \Delta^{+}_{1}}(\tte^{\alpha/2} + \tte^{-\alpha/2}).\ee
 Moreover this expression is $W$-invariant and independent of
$Z$.  Using \eqref{rat}, we apply $w$ to \eqref{asp} to get the first equality below, and replace $Z$ by $wZ$ for the second,

\by 
\tte^{\rho_{1}} \prod_{\alpha \in Z}
(1 + \tte^{-w\alpha})\sum_{\gamma \in \Gamma_Z}
K_Z(\gamma)\tte^{ -w* \gamma} &=& \prod_{\alpha \in \Delta^{+}_{1}}(\tte^{\alpha/2} + \tte^{-\alpha/2})
\nn\\
&=& \tte^{\rho_{1}} \prod_{\gb \in wZ}(1 + \tte^{-\gb})\sum_{\gamma \in \Gamma_{wZ}} K_{wZ}(\gamma)\tte^{ -\gamma}.\nn
\ey
The result follows since 

\[\tte^{\rho_{1}} \prod_{\alpha \in Z} (1 +\tte^{-w\alpha}) = \tte^{\rho_{1}} \prod_{\gb \in wZ} (1 +\tte^{-\gb}).\] 
\epf

\noi {\it Proof of Theorem \ref{hen}.}
By Lemma \ref{coq} applied to  \eqref{cod}, 
\be \ch T_w M^X(\gl)=pr_{wX}\sum_{\gc \in \Gamma_Z} K_Z(\gc)\tte^{w \circ (\lambda -\gamma)}.\ee  
Using \eqref{pig} and then Lemma \ref{pug}, we see that $T_w M^X(\gl)$ has character

\by
pr_{wX}\tte^{w\cdot \gl }\sum_{\gc \in\Gamma_{Z}}  K_Z(\gc)\tte^{-w*\gc}
&=& pr_{wX}\tte^{w\cdot \gl } \sum_{\gamma \in\Gamma_{wZ}}  K_{wZ}(\gc)\tte^{ -\gc}\nn\\
&=& \tte^{w\cdot \gl }p_{wX}.\ey
This completes the proof. 
\hfill  $\Box$
\bexa {\rm 
Here is one of the simplest cases to which Theorem \ref{hen} applies. 
Consider the following subalgebras of $\fgl(3|2),$ where the stars denote arbitrary entries

\[	
\fl=		 \begin{tabular}{|c|c||c|c|c|}\hline
  $*$&$*$&$*$&$*$&$0$\\ \hline 
$*$&$*$&$*$&$*$&$0$\\ \hline  \hline
$*$&$*$&$*$&$*$&$0$\\ \hline
$*$&$*$&$*$&$*$&$0$\\ \hline   
$0$&$0$&$0$&$0$&$*$\\ \hline 
    \end{tabular}
			\quad  \mbox{ and }\quad
 \fp =\begin{tabular}{|c|c||c|c|c|}\hline
  $*$&$*$&$*$&$*$&$*$\\ \hline 
$*$&$*$&$*$&$*$&$*$\\ \hline  \hline
$*$&$*$&$*$&$*$&$*$\\ \hline
$*$&$*$&$*$&$*$&$*$\\ \hline   
$0$&$0$&$0$&$0$&$*$\\ \hline 
    \end{tabular}  
 \]

\noi Also let $\gc_1 = \gep_2-\gd_1, \;\gc_2 = \gep_1-\gd_2, \;\ga_1 =\gep_1-\gep_2,\;\ga_2 = \gd_1-\gd_2$, 
$X=\{\gc_1, \gc_2 \}$ and $W = \Gd^+(\fl)$. 
Note that $\ga_1 +\gc_1+\ga_2 =\gc_2$, so $\ga_1^\vee \equiv \ga_2^\vee$ mod $\Z X$.
Suppose that $(\gl+\gr, \ga_1^\vee) =(\gl+\gr, \ga_2^\vee)=1, $	 and $(\gl+\gr, X)=0$.  
Then in the Jantzen sum formula for $M^X(\gl)$ the sum
$ \sum_{i>0} \ch M_i^X(\gl)$ is expressed as a positive linear combination of characters of other modules, 
one of which is $M^W(\gl-\ga_1-\ga_2):= {\Ind^\fg_\fp} \;\ttk_{\gl-\ga_1-\ga_2}$, which has character $\tte^{\gl-\ga_1-\ga_2 }p_{W}$. 
\\ \\
Now let $s$ be the reflection with respect to the root $\ga =\gd_2-\gd_3$, $Y= sX$ and $\mu = s\cdot \gl$. The Jantzen sum formula for $M^Y(\mu)$ contains a module with character $\tte^{\mu-\ga_1-\ga_2 -\ga}p_{sW}$.  Such a module cannot be obtained by parabolic induction (essentially because $\ga+\ga_2$ is not a simple root). However Theorem \ref{hen} implies that $T_s M^W(\gl-\ga_1-\ga_2)$ has this character. 

}\eexa


\begin{bibdiv}
\begin{biblist}

\bib{AL}{article}{
   author={Andersen, H. H.},
   author={Lauritzen, N.},
   title={Twisted Verma modules},
   conference={
      title={Studies in memory of Issai Schur},
      address={Chevaleret/Rehovot},
      date={2000},
   },
   book={
      series={Progr. Math.},
      volume={210},
      publisher={Birkh\"auser Boston, Boston, MA},
   },
   date={2003},
   pages={1--26},
   review={\MR{1985191 (2004d:17005)}},
}
		
\bib{AS}{article}{
   author={Andersen, Henning Haahr},
   author={Stroppel, Catharina},
   title={Twisting functors on $\scr O$},
   journal={Represent. Theory},
   volume={7},
   date={2003},
   pages={681--699},
   issn={1088-4165},
   review={\MR{2032059}},
   doi={10.1090/S1088-4165-03-00189-4},
}
				
		\bib{A}{article}{
   author={Arkhipov, Sergey},
   title={Algebraic construction of contragradient quasi-Verma modules in
   positive characteristic},
   conference={
      title={Representation theory of algebraic groups and quantum groups},
   },
   book={
      series={Adv. Stud. Pure Math.},
      volume={40},
      publisher={Math. Soc. Japan, Tokyo},
   },
   date={2004},
   pages={27--68},
   review={\MR{2074588 (2005h:17027)}},
}
		
\bib{BGG}{article}{
  author={Bern{\v{s}}te{\v{\i}}n, I. N.},
   author={Gel\cprime fand, I. M.},
   author={Gel\cprime fand, S. I.},
   title={A certain category of ${\germ g}$-modules},
   language={Russian},
   journal={Funkcional. Anal. i Prilo\v zen.},
   volume={10},
   date={1976},
   number={2},
   pages={1--8},
   issn={0374-1990},
   review={\MR{0407097}},
} 
		

\bib{CM}{article}{
   author={Coulembier, Kevin},
   author={Mazorchuk, Volodymyr},
   title={Primitive ideals, twisting functors and star actions for classical
   Lie superalgebras},
   journal={J. Reine Angew. Math.},
   volume={718},
   date={2016},
   pages={207--253},
   issn={0075-4102},
   review={\MR{3545883}},
   doi={10.1515/crelle-2014-0079},
}

\bib{H}{book}{ author={Humphreys, James E.}, title={Introduction to Lie algebras and representation theory}, note={Graduate Texts in Mathematics, Vol. 9}, publisher={Springer-Verlag}, place={New York}, date={1972}, pages={xii+169}, review={\MR{0323842 (48 \#2197)}}, }



\bib{H2}{book}{
   author={Humphreys, James E.},
   title={Representations of semisimple Lie algebras in the BGG category
   $\scr{O}$},
   series={Graduate Studies in Mathematics},
   volume={94},
   publisher={American Mathematical Society},
   place={Providence, RI},
   date={2008},
   pages={xvi+289},
   isbn={978-0-8218-4678-0},
   review={\MR{2428237}},
}

\bib{J1}{book}{ author={Jantzen, Jens Carsten}, title={Moduln mit einem h\"ochsten Gewicht}, language={German}, series={Lecture Notes in Mathematics}, volume={750}, publisher={Springer}, place={Berlin}, date={1979}, pages={ii+195}, isbn={3-540-09558-6}, review={\MR{552943 (81m:17011)}}, }


	\bib{KM}{article}{
   author={Khomenko, Oleksandr},
   author={Mazorchuk, Volodymyr},
   title={On Arkhipov's and Enright's functors},
   journal={Math. Z.},
   volume={249},
   date={2005},
   number={2},
   pages={357--386},
   issn={0025-5874},
   review={\MR{2115448}},
   doi={10.1007/s00209-004-0702-8},
}

\bib{Mz}{book}{
   author={Mazorchuk, Volodymyr},
   title={Lectures on algebraic categorification},
   series={QGM Master Class Series},
   publisher={European Mathematical Society (EMS), Z\"urich},
   date={2012},
   pages={x+119},
   isbn={978-3-03719-108-8},
   review={\MR{2918217}},
   doi={10.4171/108},
}


\bib{M0}{article}{
   author={Musson, Ian M.},
   title={On the center of the enveloping algebra of a classical simple Lie
   superalgebra},
   journal={J. Algebra},
   volume={193},
   date={1997},
   number={1},
   pages={75--101},
   issn={0021-8693},
   review={\MR{1456569 (98k:17012)}},
   doi={10.1006/jabr.1996.7000},
}




\bib{M}{book}{author={Musson, Ian M.}, title={Lie Superalgebras and Enveloping Algebras},
   series={Graduate Studies in Mathematics},
   volume={131},
publisher={American Mathematical Society}, place={Providence, RI}, date ={2012}}



		



\bib{M1001}{article}{author={Musson, Ian M.}, title={\v Sapovalov elements and the Jantzen sum formula for contragredient Lie superalgebras.}, 
journal={arXiv 1710.10528.}}

\bib{P}{article}{
   author={Pedersen, Dennis Hasselstr\o m},
   title={Twisting functors for quantum group modules},
   journal={J. Algebra},
   volume={447},
   date={2016},
   pages={580--623},
   issn={0021-8693},
   review={\MR{3427652}},
   doi={10.1016/j.jalgebra.2015.09.046},
}
		


\end{biblist}

\end{bibdiv}
 \end{document}